\documentclass[12pt]{amsart}
\usepackage{amsmath}
\usepackage{amssymb}
\usepackage{amscd}
\usepackage{verbatim}
\newcommand{\be}{\begin{equation}}
\newcommand{\bs}{\begin{sub}}
\newcommand{\es}{\end{sub}}
\newcommand{\bsn}{\begin{subn}}
\newcommand{\esn}{\end{subn}}
\newcommand{\bea}{\begin{eqnarray}}
\newcommand{\eea}{\end{eqnarray}}
\newcommand{\BA}[1]{\begin{array}{#1}}
\newcommand{\EA}{\end{array}}
\newcommand{\Real}{\mathbb{R}}

\newcommand{\ZZ}{\mathbb{Z}}
\newcommand{\CC}{\mathbb{C}}
\newcommand{\RR}{\mathbb{R}}
\newcommand{\g}{\gamma}
\newcommand{\G}{\Gamma}

\newlength{\wex}  \newlength{\hex}
\def\ga{\alpha}     \def\gb{\beta}

            \def\gl{\lambda}
\def\gm{\mu}

\newcommand{\pf}{\noindent \mbox{{\bf Proof}: }}
\def\squarebox#1{\hbox to #1{\hfill\vbox to #1{\vfill}}}

\newcommand{\beq}{\begin{equation}}
\newcommand{\eeq}{\end{equation}}
\newcommand{\beqa}{\begin{eqnarray}}
\newcommand{\eeqa}{\end{eqnarray}}
\newcommand{\beqanl}{\begin{eqnarray*}}
\newcommand{\eeqanl}{\end{eqnarray*}}
\begin{document}
\renewcommand{\theequation}{\thesection.\arabic{equation}}
\newcommand{\mysection}[1]{\section{#1}\setcounter{equation}{0}}

\def\stackunder#1#2{\mathrel{\mathop{#2}\limits_{#1}}}
\newtheorem{theorem}{Theorem}
\newtheorem{lemma}[theorem]{Lemma}
\newtheorem{definition}[theorem]{Definition}
\newtheorem{corollary}[theorem]{Corollary}
\newtheorem{conjecture}[theorem]{Conjecture}
\newtheorem{remark}[theorem]{Remark}
\def\stackunder#1#2{\mathrel{\mathop{#2}\limits_{#1}}}

\author{Peter Kuchment}
\address{Department of Mathematics\\ Texas A\&M University\\
College Station, TX 77843-3368, USA}
\email{kuchment@math.tamu.edu}
\title[Integral representations]%
{Integral representations of solutions of periodic elliptic
equations} \dedicatory{Dedicated to Stas Molchanov on the occasion
of his 65th birthday}
\date{}
\begin{abstract}
The paper discusses relations between the structure of the complex
Fermi surface below the spectrum of a second order periodic
elliptic equation and integral representations of certain classes
of its solutions. These integral representations are analogs of
those previously obtained by S.~Agmon, S.~Helgason, and other
authors for solutions of the Helmholtz equation (i.e., for
generalized eigenfunctions of Laplace operator). In a previous
joint work with Y.~Pinchover we described all solutions that can
be represented as integrals of positive Bloch solutions over the
imaginary Fermi surface, with a hyperfunction as a ``measure''.
Here we characterize the class of solutions such that the
corresponding hyperfunction is a distribution on the Fermi
surface.

\end{abstract}

\subjclass[2000]{Primary: 35B05, 35C15, 58J15; Secondary: 35J15,
35P05, 58J50.}

\keywords{Elliptic operator, spectrum, Floquet theory, integral
representation, periodic operator} \maketitle

\mysection{Introduction}\label{Sec1}

This paper is devoted to integral representations of solutions of
second order elliptic periodic differential equations. These
representations are analogs of those for solutions of the
Helmholtz equation in $\Real^n$
\begin{equation}
-\Delta u-k^2u=0 \quad \mbox{ in } \Real^n\, ,  \label{Helm}
\end{equation}
where $k\in \CC^{*}:=\CC\setminus\{0\}$. Such representations have
been obtained by S.~Helgason \cite{H1,H2}, M.~Hashizume et al.
\cite{HKMO}, M.~Morimoto \cite{Mo}, and S.~Agmon
\cite{Agmon_helm1,Agmon_helm2}. In these results, solutions were
expanded into exponential ones
\begin{equation}\label{E:exponent}
e_{\,\xi }(x):=\exp (\mathrm{i}\xi \cdot x).
\end{equation}
Here
\begin{equation}\label{E:sphere}
  \xi \in S:=\left\{ \xi \in \CC^{\, n}|\,\xi =k\omega ,\,\omega \in
S^{n-1}\subset \mathbb{R}^n\right\},
\end{equation}
and $\xi \cdot x=\sum_{j=1}^n \xi_jx_j.$ These expansions can be
written as
\begin{equation}\label{E:const_coef_expans}
u(x)=\int_S e_{\,\xi }(x)\,\mathrm{d}\phi (\xi )=<\phi (\xi
),e_{\,\xi }(x)>,
\end{equation}
where $\phi (\xi )$ is a functional on the sphere $S$. In
particular, it was understood what classes of solutions correspond
to different classes of functionals (e.g., hyperfunctions,
distributions, measures) \cite{Agmon_helm1,Agmon_helm2}.

Such representations are related to the L.~Ehrenpreis' {\em
fundamental principle} \cite{E,Pa1} for constant coefficient
operators, which in the particular case of (\ref{Helm}) claims
that any solution of (\ref{Helm}) can be represented as an
integral with respect to the parameter $\xi$ of the exponential
solutions
\begin{equation}\label{E:fund_princ}
e_{\,\xi }(x):=\exp (\mathrm{i}\xi \cdot x)\qquad \xi \in \Sigma.
\end{equation}
Here
$$
\Sigma :=\left\{ \xi \in \CC^{\, n}|\,\xi ^2=k^2\right\} ,
$$
is the characteristic variety of the operator in the left hand
side of (\ref{Helm}) (see the details and more precise formulation
in \cite{E,Pa1}). The representation (\ref{E:fund_princ}) is
highly non-unique, due to existence of functionals orthogonal to
analytic functions on $\Sigma$. On the other hand, $\Sigma$ is an
analytic subset of $\CC^{\, n}$, uniquely determined for $k\neq 0$
by its spherical subset $S$. Thus, one can expect the possibility
of a unique representation like the one in
(\ref{E:const_coef_expans}). It is crucial here that $\Sigma $ is
irreducible and that $S$ is sufficiently massive, so $S$
determines $\Sigma $ uniquely (otherwise it would not be possible
to obtain the representation of all solutions using only $\xi \in
S$). Moreover, $S$ is a rather simple analytic manifold. This
enables one to obtain rather explicit descriptions of the needed
spaces of test functions and functionals.

In this paper, we consider the case of second order periodic
elliptic equations (see the exact description of the class of
equations in the next section). For such (and more general)
periodic equations, an analog of the ``fundamental principle'' was
obtained in \cite{Ku,Pa2} for solutions with some growth
restrictions. Here, instead of exponential solutions one needs to
use the so called {\em Floquet-Bloch} solutions. The analog of the
characteristic manifold $\Sigma $ is the Fermi surface $F$  (see
\cite{AM,Ku,RS} and definitions \ref{deffloq} and \ref{defFermi}
below for these notions). However, finding an analog of
(\ref{E:const_coef_expans}) for the periodic case is far from
being straightforward. In particular, one wonders what should be
the analog of the sphere $S$. It is natural to expect that when
zero belongs to the spectrum of the operator, one might try to use
the {\em real part} of the Fermi surface, while if zero does not
belong to the spectrum, the {\em imaginary part} might be
appropriate. Due to the complicated structure of the Fermi
surface, this idea is not easy to implement. As it was shown in
the joint work \cite{KuPin} with Y.~Pinchover, for second order
elliptic equations with {\em positive} generalized principal
eigenvalue $\Lambda_0$ (see (\ref{genev})), an appropriate variety
is provided by the analysis of the cone of positive solutions done
by S. Agmon and by V. Lin and Y. Pinchover
\cite{Agmon_positive,LP,Ku}. No results of this type are known so
far above $\Lambda_0$. Another difficulty is in proving
irreducibility of the Fermi surface $F$, which happens to be a
very hard problem (it also arises in many other spectral
considerations \cite{BKT,GKT,KT,Ku_CONM,KV1,KV2,KV3}).
Fortunately, as the reader can see from this paper and from
\cite{KuPin}, by appropriately restricting the growth of the
solutions, we manage to work near a single irreducible component
of $F$ and hence avoid proving the irreducibility of $F$.
Consequently, we prove a representation theorem (Theorem
\ref{T:distribution}) that characterizes classes of the solutions
that have integral expansion analogous to
(\ref{E:const_coef_expans}) into {\em positive} Bloch solutions
with hyperfunctions or distributions as ``measures''. The
hyperfunction case was investigated in \cite{KuPin} and is
presented here without a proof. The distribution result that we
could not obtain in \cite{KuPin} is new and is proved in the
present paper. In order to prove it, additional analytic
techniques need to be involved, in particular results on peak sets
in $A^\infty$ functional algebras in several complex variables
\cite{ChaCho,ChaCho2}. We are grateful to A.~Tumanov for pointing
us to the relevant literature.

The proofs of the results of this paper are based upon the
techniques of the Floquet theory developed in \cite{Ku} (the
reader can find all necessary preliminary information in the next
section). The methods that were used to prove the ``fundamental
principle'' \cite{E,Pa2} provide a crucial technical approach. In
particular, solutions of the equation are treated in the dual
sense, i.e., as functionals on appropriate function spaces that
are orthogonal to the range of the dual operator.

The outline of the paper is as follows. The next section
introduces necessary notations and preliminary results from the
Floquet theory and the theory of positive solutions of periodic
elliptic equations. It mostly (but not entirely) repeats the
corresponding sections from \cite{KuPin,KuPin2} and is included
for the reader's convenience. Section \ref{Sec3} contains the
proof of the integral representation (Theorem
\ref{T:distribution}) that describes the sets of solutions
allowing integral representations with distributional and
hyperfunction ``measures''. The last Section contains
acknowledgments.

\mysection{Notations and preliminary results}\label{notations}

Due to the nature of this section, most of it repeats some parts
of \cite{Ku,KuPin,KuPin2}. We regret the necessity of doing this,
but otherwise reading the rest of the paper would probably become
impossible without constant referring to \cite{KuPin}.

In this paper we consider second order elliptic operators on
$\mathbb{R}^n$ with {\it real} periodic coefficients of the form
\begin{equation}\label{operator}
P(x,\partial)=-\sum_{i,j=1}^n a_{ij}(x)\partial _{i}\partial
_{j}+\sum_{i=1}^n b_i(x)\partial _i+c(x).
\end{equation}
It is assumed that the uniform ellipticity condition
$$
\sum_{i,j=1}^n a_{ij}(x)\zeta _i\zeta _j\geq
a\sum_{i=1}^n\zeta_i^2
$$
is satisfied for all $x, \zeta \in \Real^n$, where $a$ is a
positive constant. In the notation $P(x,D)$ we used the standard
convention $D=-\mathrm{i}{\partial}/{\partial x}$.

We will assume sufficient smoothness of the coefficients, namely
that $a_{ij}\in C^2(\RR^n)$, $b_i\in C^1(\RR^n)$ and $c\in
C(\RR^n)$. In fact, it is sufficient to assume that both the
operator $P$ and its dual (the formal adjoint) $P^*$ have
H\"{o}lder continuous coefficients\footnote{See \cite[Section
6]{KuPin}. We only need that both operators $P$ and $P^{*}$ define
Fredholm mappings between the Sobolev space $H^{2}(\mathbb{T}^n)$
and $L_2(\mathbb{T}^n)$, where $\mathbb{T}^n=\Real^{n} /
\Gamma$.}. Here the duality is provided by the bilinear (rather
than the sesquilinear) form
$$\label{E:duality}
<g,f>:=\int\limits_{\Real^n}f(x)g(x)\,\mathrm{d}x.
$$
So, the dual operator $P^{*}$ has similar properties to the ones
of $P$.

The coefficients of $P$ are assumed to be periodic with respect to
a lattice $\Gamma $ in $\Real^n$. In what follows, the particular
choice of the lattice is irrelevant and can always be reduced by
change of variables to the case $\Gamma =\ZZ^n$, which we will
assume from now on. We will always use the word ``periodic" in the
meaning of ``$\Gamma$-periodic".

We now introduce some standard notions and results from Floquet
theory of periodic PDEs \cite{AM,Ea,Ku,KuPin,KuPin2,RS}.

We denote by $K=[0,1]^n\,$ the standard fundamental domain (the
{\it Wigner-Seitz cell}) of the lattice $\Gamma =\ZZ^n$, and by
$B=[-\pi ,\pi ]^n$ the {\it first Brillouin zone}, which is a
fundamental domain of the reciprocal (dual) lattice $\Gamma
^{*}=\left( 2\pi \ZZ\right) ^n$. We naturally identify
$\Gamma$-periodic functions with functions on ${\mathbb
T}^{n}=\Real^{n} / \Gamma$.
\begin{definition}\label{deffloq}
{\em A function $u(x)$ representable as a finite sum of the form
\begin{equation}\label{floquetsol}
u(x)=\mathrm{e}^{\mathrm{i}k\cdot x}\left( \sum\limits_{j=(j_1,\ldots ,j_n)\in \ZZ%
_{+}^n}x^jp_j(x)\right)
\end{equation}
with nonzero $\G$-periodic functions $p_j(x)$ is called a {\em
Floquet function with a quasimomentum $k\in \CC^n$}. Here
 $x^j =x_1^{j _1}x_2^{j_2}\ldots x_n^{j _n}$. The
maximum value of $|j|=\sum\limits_{l=1}^nj_l$ in the
representation (\ref{floquetsol}) is said to be the {\em order} of
the Floquet function. Floquet functions of zero order are called
{\em Bloch functions}.}
\end{definition}

The set introduced below plays in the periodic case the role of
the characteristic variety $\Sigma $ for constant coefficient
operators.
\begin{definition} \label{defFermi}
{\em The (complex) {\em Fermi surface} $F_P$ of the operator $P$
(at the zero energy level) consists of all vectors $k\in \CC^{\,
n}$ (called {\em quasimomenta}) such that the equation $Pu=0$ has
a nonzero Bloch solution $u(x)=\mathrm{e}^{\mathrm{i}k\cdot
x}p(x)$, where
 $p(x)$ is a $\Gamma $-periodic function.}
\end{definition}

Introducing a spectral parameter $\lambda$, one arrives at the
notion of the {\it Bloch variety}:
\begin{definition} \label{defBloch}
{\em The (complex) {\em Bloch variety} $B_P$ of the operator $P$
consists of all pairs $(k,\lambda) \in \CC^{\, n+1}$ such that the
equation $Pu=\lambda \, u$ has a nonzero Bloch solution
$u(x)=\mathrm{e}^{\mathrm{i}k\cdot x}p(x)$ with the quasimomentum
$k$.}
\end{definition}

The Bloch variety $B_P$ can be treated as the graph of a
multivalued function $\lambda(k)$ (so called {\it dispersion
relation}) that assigns to any quasimomentum $k$ the spectrum of
the operator $P(x,D+k)$ on the torus $\mathbb{T}^n$. Since for
operators of the type (\ref{operator}), these spectra are known to
be discrete (as in particular the discussion below will show), we
can single out continuous branches $\lambda_j$ of $\lambda(k)$.
These branches are usually called the {\it band functions} (see
\cite{AM, RS, Ku}). The Fermi surfaces now become the level sets
of the dispersion relation.
\begin{lemma}[{\cite[Theorems 3.1.5, 3.1.7 and 4.4.2]{Ku}}]
\label{Fermi}
\begin{enumerate}
\item The Fermi and Bloch varieties are the sets of all zeros of
entire functions of a finite order in $\CC^{n}$ and $\CC^{n+1}$,
respectively.
\item A quasimomentum $k$ belongs to $F_{P^*}$ if
and only if $-k \in F_{P}$. Analogously, $(k,\lambda)\in B_{P^*}$
if and only if $(-k,\lambda)\in B_{P}$. In other words, the
dispersion relations $\lambda(k)$ and $\lambda^*(k)$ for the
operators $P$ and $P^*$ are related as follows:
\begin{equation}
\lambda^*(k)=\lambda(-k).
\label{dual_disp}
\end{equation}
\end{enumerate}
\end{lemma}

The Fermi surface $F_P$ is periodic with respect to the reciprocal
lattice $\Gamma ^{*}=(2\pi \ZZ)^n$. It is often convenient to
factor out this periodicity by considering the (analytic)
exponential mapping $\rho :\CC^{n}\rightarrow (\CC^{*})^n$, where
$$
z=\rho (k)=\rho (k_1,\ldots,k_n)=(\exp \mathrm{i}k_1,\ldots,\exp
\mathrm{i}k_n).
$$
This mapping can be identified with the quotient map
$\CC^{n}\rightarrow \CC^{n}/\Gamma ^{*}$. We also introduce the
complex torus
\begin{equation}
T:=\rho (\Real^n)=\left\{ z\in \CC^{n}|\,\left| z_j\right|
=1,\,j=1,2,\dots,n\right\} .  \label{torus}
\end{equation}

\begin{definition}
\label{def_Floquet} {\em We call the image $\Phi_P:=\rho (F_P)$
under the mapping $\rho$ of the Fermi surface $F_P$ the {\em
Floquet surface} of the operator $P$.}
\end{definition}

In the Floquet theory for PDEs, this Floquet surface is the set of
all Floquet multipliers of Floquet-Bloch solutions of the equation
$Pu=0$.

The following analog ${\mathcal U}$ of the Fourier transform (see
\cite[Section 2.2]{Ku}, \cite{RS}), which we will call the {\it
Floquet transform}\footnote{It is sometimes also called the {\it
Gelfand transform}, due to Gelfand's work \cite{Gelf}.}, is the
main tool in the Floquet theory:
\begin{equation}
f(x)\rightarrow {\mathcal U}f(z,x):=\sum_{\gamma \in \Gamma
}f(x-\gamma )z^\gamma \qquad z\in (\CC^{*})^n.  \label{Floquet}
\end{equation}

It is often convenient to use for the Floquet transform ${\mathcal
U}$ the quasimomentum coordinate $k$ instead of the multiplier
$z=\rho (k)$.

For a point $z\in (\CC^{*})^n$, we denote by $E_{m,z}$ the closed
subspace of the Sobolev space $H^m(K)$ formed by the restrictions
of functions $v\in H_{\mathrm{loc}}^m(\Real^n)$ that satisfy the
Floquet condition $v(x+\gamma )=z^\gamma v(x)$ for any $\gamma \in
\Gamma $. One can show \cite[Theorem 2.2.1]{Ku} that
\begin{equation}
{\mathcal E}_m:=\mathrel{\mathop{\cup }\limits_{z\in (\CC^{*})^n}}%
E_{m,z}  \label{bundle}
\end{equation}
forms a holomorphic subbundle of the trivial bundle
$(\CC^{*})^n\times H^m(K)$. As any infinite dimensional analytic
Hilbert bundle over a Stein domain, it is trivializable (see
\cite[Chapter 1]{Ku} and Lemma \ref{Plancherel} below). One can
notice that for $m=0$ the bundle ${\mathcal E}_0$ coincides with
the whole $(\CC^{*})^n\times L^2(K)$.

We collect now several statements from Theorem {\small{\sf
XIII}}.97 in \cite{RS} and Theorems  1.3.2, 1.3.3,  1.5.23 and
2.2.2 in \cite{Ku}:
\begin{lemma}\label{Plancherel}
\begin{enumerate}
\item As any infinite dimensional analytic Hilbert bundle over a
Stein domain, the bundle ${\mathcal E}_m$ is analytically trivial.

\item  For any nonnegative integer $m$, the operator
$$ {\mathcal
U}:H^{m}(\Real^n)\rightarrow L^2(T,{\mathcal E}_m)
$$
is an isometric isomorphism, where $L^2(T,{\mathcal E}_m)$ denotes
the space of square integrable sections over the complex torus $T$
of the bundle ${\mathcal E}_m$, equipped with the natural topology
of a Hilbert space.

\item Let the space
$$ \Theta ^m:=\left\{f\in H_{\mathrm{loc}}^m(R^n)|\,
\underset{\gamma \in \Gamma }{\sup}\{||f||_{H^m(K+\gamma
)}\mathrm{e}^{b|\gamma|}\}<\infty \,,\; \forall b>0\right\}$$ be
equipped with the natural Fr\'{e}chet topology. Then
$$
{\mathcal U}:\Theta^{m}\rightarrow \Gamma((\CC^{*})^{n},{\mathcal
E}_m)
$$
is a topological isomorphism, where
$\Gamma((\CC^{*})^{n},{\mathcal E}_m)$ is the space of all
analytic sections over $(\CC^{*})^{n}$ of the bundle ${\mathcal
E}_m$, equipped with the topology of uniform convergence on
compacta.

\item Under the transform ${\mathcal U}$, the operator
$$
P: H^{2}(\Real^n) \rightarrow L^{2}(\Real^n)
$$
becomes the operator of multiplication by a holomorphic Fredholm
morphism $P(z)$ between the fiber bundles ${\mathcal E}_2$ and
${\mathcal E}_0$. Here $P(z)$ acts on the fiber of ${\mathcal
E}_m$ over the point $z \in T$ as the restriction to this fiber of
the operator $P$ acting between $H^2(K)$ and $L^2(K)$.
\end{enumerate}
\end{lemma}

Let us now mention another common way of looking at $P(z)$. If
$z=\exp \mathrm{i}k$, then commuting with the exponent $\exp
\mathrm{i}k \cdot x$ one reduce the bundle ${\mathcal E}_m$ to the
trivial one with the fiber $H^{m}({\mathbb T}^n)$, where as before
${\mathbb T}^n=\Real^n / \Gamma$. On the other hand, the operator
$P(z)$ takes the form $P(x,D+k)$ acting between Sobolev spaces on
the torus ${\mathbb T}^n$. In other words, the options are either
dealing with the restriction of a fixed operator to an
analytically ``rotating'' subspace, or with a polynomial family of
operators between fixed spaces.

We will need to see how the structure of the Floquet solutions
(see Definition \ref{deffloq}), and in general, the structure of
functions of Floquet type (\ref{floquetsol}) reacts to the Floquet
transform. For instance, in the constant coefficient case, where
the role of the Floquet solutions is played by the exponential
polynomials
$$ \mathrm{e}^{\mathrm{i}k\cdot x}\sum\limits_{\left| j\right| \leq N}p_jx^j,$$ such
functions are Fourier transformed into distributions supported at
the point
 $\left(-k\right)$. The next statement shows that
under the Floquet transform, each Floquet type function of the
form (\ref{floquetsol}) corresponds, in a similar way, to a
(vector valued) distribution supported at the quasimomentum
$\left( -k\right) $.

Every Floquet type function $u$ (see (\ref{floquetsol})), being of
exponential growth, determines a (continuous linear) functional on
the space $\Theta ^0 $. If it satisfies the equation $Pu=0$ for a
periodic elliptic operator of order $m$, then as such a functional
it is orthogonal to the range of the dual operator $P^{*}:\Theta
^m\rightarrow \Theta ^0$. According to Lemma \ref{Plancherel},
after the Floquet transform any such functional becomes a
functional on $\Gamma \left( \left( \CC^{*}\right) ^n,{\mathcal
E}_0\right)$, which is orthogonal to the range of the Fredholm
morphism $P^{*}(z):{\mathcal E}_m\rightarrow {\mathcal E}_0$
generated by the dual operator $P^{*}:\Theta ^m\rightarrow \Theta
^0$. The following auxiliary result describes all such
functionals.
\begin{lemma}[{\cite[Lemma 8]{KuPin}}]\label{L:Floq_struct}
A continuous linear functional $u$ on $\Theta ^0$ is generated by
a function of the Floquet form (\ref{floquetsol}) with a
quasimomentum $k$ if and only if after the Floquet transform it
corresponds to a functional on $\Gamma \left( \left(
\CC^{*}\right) ^n,{\mathcal E}_0\right) $ which is a distribution
$\phi $ that is supported at the point $\nu =\exp (-\mathrm{i}k)$,
i.e. has the form $$ \left\langle \phi ,f\right\rangle
=\sum_{\left| j\right| \leq N}\left\langle q_j,\left.
\frac{\partial ^{\left| j\right| }f}{\partial z^j}\right| _{\,\nu
}\right\rangle  \qquad f\in \Gamma \left( \left( \CC^{*}\right)
^n,{\mathcal E}_0\right),
$$ where $q_j\in L^2(K)$. The orders $N$ of the Floquet function
(\ref{floquetsol}) and of the corresponding distribution $\phi $
are the same.
\end{lemma}
Everything discussed so far applies to essentially any elliptic
periodic scalar or matrix operators of any order, not necessarily
to the ones of the form (\ref{operator}) (see
\cite{Ku,KuPin,KuPin2}). However, there is a special construction
that applies only to operators (\ref{operator}) and which will
play a crucial role in our considerations. Its properties were
studied in detail in \cite{Agmon_positive,LP,Pr}.

Consider the function $\Lambda (\xi ):\RR^n\rightarrow \RR$
defined by the condition that the equation
$$
Pu=\Lambda (\xi )u
$$
has a {\em positive} Bloch solution of the form
\begin{equation}\label{positiveBloch}
u_{\,\xi }(x)=\mathrm{e}^{\xi \cdot x}p_{\,\xi }(x),
\end{equation}
where $p_{\,\xi }(x)$ is $\Gamma $-periodic.
\begin{lemma}[{[Lemma 12]\cite{KuPin}}]
\label{Lambda-lemma}
\begin{enumerate}
\item  The value $\Lambda (\xi )$ is uniquely determined for any $\xi \in
\Real^n$.

\item  The function $\Lambda (\xi )$ is bounded from above, strictly
concave, analytic, and has a nonzero gradient at all points except
at its maximum point.

\item  Consider the operator
$$
P(\xi )=\mathrm{e}^{-\xi\cdot x}P\mathrm{e}^{\xi\cdot
x}=P(x,D-\mathrm{i}\xi)
$$
on the torus $\mathbb{T}^n$.  Then $\Lambda (\xi )$ is the
principal eigenvalue of $P(\xi )$ with  a positive eigenfunction
$p_{\,\xi }$. Moreover, $\Lambda (\xi )$ is algebraically simple.

\item The Hessian of $\Lambda (\xi )$ is nondegenerate at all points.

\end{enumerate}
\end{lemma}

One should note that since the function $\Lambda(\xi)$ is analytic,
it is actually defined in a neighborhood of $\Real^n$ in $\CC^n$.
This remark will be used in what follows.

Let us denote
\begin{equation}
\Lambda_{0} :=\max_{\xi \in \Real^n}\Lambda (\xi ).
\label{Lambda}
\end{equation}
It follows from \cite{Agmon_positive,LP} that an alternative
definition of $\Lambda_{0}$ can be
\begin{equation}\label{genev}
\Lambda_0= \sup\{\gl \in \Real \; |\; \exists u>0 \mbox{ such that
} (P-\gl)u= 0 \mbox{ in } \Real^n\},
\end{equation}
and that in the self-adjoint case $\Lambda_{0}$ coincides with the
bottom of the spectrum of the operator $P$. The common name for
$\Lambda_{0}$ is {\em the generalized principal eigenvalue} of the
operator $P$ in $\Real^n$.

In our main result, we will need to assume that $\Lambda_{0}$ is
strictly positive. In the self-adjoint case such an assumption has
a clear spectral interpretation: the bottom of the spectrum is
strictly positive. In the next lemma, we provide some known
conditions for the nonnegativity or positivity of $\Lambda_{0}$
for not necessarily self-adjoint operators of the form
(\ref{operator}).
\begin{lemma}[{\cite[Lemma 13]{KuPin}}]
\label{Lambdazero} Consider an operator $P$ of the form
(\ref{operator})
\begin{enumerate}

\item $\Lambda _{0}\geq 0$ if and only if the operator $P$ admits
a positive (super)solution. This condition is satisfied in
particular when $c(x)\geq 0$.

\item $\Lambda _{0}\geq 0$ if and only if the operator $P$ admits
a positive solution of the form (\ref{positiveBloch}).

\item $\Lambda _{0}=0$ if and only if the equation $Pu=0$ admits
exactly one normalized positive solution in $\Real^n$.

\item If $c(x)=0$, then $\Lambda_{0} =0$ if and only if
$\int\limits_{{\mathbb T}^n}b(x)\psi (x)\,\mathrm{d}x=0$, where
$\psi $ is the principal eigenfunction of $P^{*}$ on ${\mathbb
T}^n$ (with principal eigenvalue zero). In particular, divergence
form operators satisfy this condition.

\item Let $\xi\in \Real^n$, and assume that  $u_\xi(x)=\mathrm{e}^{\xi
\cdot x}p_\xi (x)$ and $u^*_{-\xi}$ are positive Bloch solutions
of the equations $Pu=0$ and $P^*u=0$, respectively. Denote by
$\psi$ the periodic function $u_\xi u^*_{-\xi}\,$. Consider the
function
$$
\tilde{b}_i(x):= b_i(x)-2\sum_{j=1}^n a_{ij}(x)\{\xi_j
+[p_\xi(x)]^{-1}\partial_{j}p_\xi(x)\},
$$
and denote
$$
\gamma =(\gamma _1,\ldots ,\gamma
_n):=(\int\limits_{\mathbb{T}^n}\tilde{b}_1(x)\psi
(x)\,\mathrm{d}x,\ldots
,\int\limits_{\mathbb{T}^n}\tilde{b}_n(x)\psi (x)\,\mathrm{d}x).
$$
Then $\Lambda_{0} =0$ if and only if $\gamma =0$.

\end{enumerate}
\end{lemma}
Let us discuss also  some additional properties that will play an
important role in the sequel. Assume that $\Lambda_{0}> 0$. Then
Lemma \ref{Lambda-lemma} implies that the zero level set
\begin{equation}
\Xi :=\left\{ \xi \in \Real^n|\;\Lambda (\xi )=0\right\}
\label{ksi}
\end{equation}
is a strictly convex compact analytic surface in $\Real^n$ of
dimension $n-1$. The manifold $\Xi$ consists of all $\xi \in
\Real^n$ such that the equation $Pu=0$ admits a positive Bloch
solution $u_{\,\xi}(x)=\mathrm{e}^{\xi\cdot x}p_{\,\xi }(x)$.
Moreover, the set of all such positive Bloch solutions is the set
of all {\em minimal positive solutions} of the equation $Pu=0$ in
$\mathbb{R}^n$ \cite{Agmon_positive,LP}\footnote{It is also
established that a function $u$ is a positive solution of the
equation $Pu=0$ in $\mathbb{R}^n$ if and only if there exists a
positive finite measure $\gm$ on $\Xi$ such that
$$
u(x)=\int_\Xi u_{\,\xi }(x)\,\mathrm{d}\gm(\xi ).
$$
}. We denote by $G$ the convex hull  of $\Xi$, and by
$\stackrel{\circ}{G}$ its interior ($\stackrel{\circ}{G}$ is
nonempty if and only if $\Lambda_{0}>0$).
\begin{lemma}[{\cite[Lemma 14]{KuPin}}]
\label{analytic} Suppose that $\Lambda_{0}>0$. There exists a
neighborhood $W$ of $G$ in $\CC^n$ and an analytic function
$$
W\ni\xi \mapsto p_{\xi}(\cdot) \in H^2(\mathbb{T}^n )
$$
such that for any $\xi \in W$ the function of $x$
$$
u_\xi (x)=\exp(\xi \cdot x)p_{\xi}(x)
$$
is a nonzero Bloch solution of the equation $Pu=\Lambda (\xi)u$
with a quasimomentum $-\mathrm{i} \xi$. Moreover, one can choose
the function $p$ in such a way that it is positive for all $\xi\in
\Xi$.
\end{lemma}

Comparing $\Xi $ with the Fermi surface $F_P$, one sees that
$$
-\mathrm{i}\Xi \subset F_P.
$$
The next result specifies further the relation between these two
varieties:
\begin{lemma}[{\cite[Lemma 15]{KuPin}}]
\label{Fermi-lemma} Let $\Lambda_{0} \geq0$. Then
\begin{enumerate}
\item  The intersection of the complex Fermi surface $F_P$ with
the tube
\begin{equation}
{\mathcal T}:=\left\{ k\in \CC^{n}|\;\mathrm{Im} \, k=(\mathrm{Im}
\, k_1,\dots ,\mathrm{Im} \, k_n)\in -G\right\} \label{tube}
\end{equation}
coincides with the union of the surface $-\mathrm{i}\Xi $ with its
translations by the vectors of the reciprocal lattice $\Gamma
^{*}$, i.e. consists of vectors $k=-\mathrm{i}\xi +\gamma $ where
$\xi \in \Xi $ and $\gamma \in \Gamma ^{*}$. Moreover, up to a
multiplicative constant, any nonzero Bloch solution with a
quasimomentum in the above intersection is a positive Bloch
solution.

\item  If $\Lambda_{0} >0$, then the intersection of $F_P$ with a
sufficiently small neighborhood of $-\mathrm{i}\Xi $ is a (smooth)
analytic manifold that coincides with the set of zeros of the
function $\Lambda(\mathrm{i}k)$.
\end{enumerate}
\end{lemma}

Analogously to the Floquet surface $\Phi=\Phi_P$, we define the
surface
\begin{equation}
\Psi :=\rho (-\mathrm{i}\Xi )=\left\{ z\,|\; z=(\exp \xi_1,\ldots,
\exp \xi_n),\;\;\xi\in \Xi\right\}, \label{Psi}
\end{equation}
and the tubular domain
\begin{equation}\label{tubulard}
V:=\rho({\mathcal T}),
\end{equation}
where ${\mathcal T}$ was defined in (\ref{tube}). The results of
lemmas \ref{analytic} and \ref{Fermi-lemma} can be restated in
terms of these new objects:

\begin{lemma}
\label{Lambda-lemma2} Let $\Lambda_0\geq 0$. Then
\begin{enumerate}
\item  $\Phi \cap V=\Psi $.

If $\Lambda_{0} >0$, then

\item  The intersection of $\Phi $ with a sufficiently small neighborhood of
$\Psi $ is a (smooth) connected analytic manifold.

\item  The intersections of $\Phi $ with neighborhoods of the tube $V$ form
a basis of neighborhoods of $\Psi $ in $\Phi $.

\item  For a sufficiently small neighborhood $\Phi _{\,\varepsilon }$ of
$\Psi $ in $\Phi $ there exists an analytic function
$p:\Phi_{\,\varepsilon }\rightarrow H^2(\mathbb{T}^n)$ such that
for any $z\in \Phi _{\,\varepsilon }$ the function of $x$
$$
u_z(x)=z^xp(z,x)
$$
is a nonzero Bloch solution of the equation $Pu=0$.
\end{enumerate}
\end{lemma}

\mysection{Representation of solutions by hyperfunctions and
distributions}\label{Sec3}

The main result of this paper (Theorem \ref{T:distribution} below)
is analogous to the results of \cite{Agmon_helm1,Agmon_helm2} that
characterize the classes of solutions of the Helmholtz equation
that can be represented by means of distributions or
hyperfunctions on $S$ (see also the introduction to our paper). In
order to state it, we need to introduce a new object. Let us
denote by $h(\omega )$,  $\omega \in S^{n-1}$ the indicator
function of the convex domain $G$ introduced in the previous
section. Namely,
\begin{equation}
h(\omega ):=\mathrel{\mathop{\sup }\limits_{\xi\in G}}(\omega
\cdot \xi), \label{indicator}
\end{equation}
where $\omega \cdot \xi=\sum_{j=1}^n \omega _j\xi_j$ is the inner
product in $\Real^n$. The next main theorem will be stated in
terms of this function.

\begin{theorem} \label{T:distribution}
Suppose that $\Lambda_{0} >0$.
\begin{enumerate}
\item Let  $u$ be a solution of the equation $Pu=0$ in $\Real^n$
satisfying for some $N$ the estimate
\begin{equation}\label{E:distr_estimate}
\left| u(x)\right| \leq C(1+|x|)^{N}\mathrm{e}^{h(x/\left|
x\right| ) \left| x\right|}.
\end{equation}
Then $u$ can be represented as
\begin{equation}
u(x)=<\mu (\xi ),u_{\,\xi }(x)>,  \label{E:represent}
\end{equation}
where $u_{\,\xi} $ is the analytic positive Bloch solution
corresponding to $\xi \in \Xi$ (see Lemma \ref{analytic}), and $
\mu (\xi )$ is a distribution on $\Xi $. The converse statement is
also true: for any distribution $\mu $ on $\Xi $, the function
$u(x)$ in (\ref{E:represent}) is a solution of the equation $Pu=0$
in $\Real^n$ which satisfies for some $N$ the growth condition
(\ref{E:distr_estimate}). \item Let  $u$ be a solution of the
equation $Pu=0$ in $\Real^n$ satisfying for any $\varepsilon >0$
the estimate
\begin{equation}\label{E:hyper_estimate}
\left| u(x)\right| \leq C_{\,\varepsilon }\exp \left[ \left(
h(x/\left| x\right| )+\varepsilon \right) \left| x\right| \right],
\end{equation}
where $C_{\,\varepsilon }$ is a constant depending only on
$\varepsilon$ and $u$. Then $u$ can be represented as in
(\ref{E:represent}) with $ \mu (\xi )$ being a hyperfunction
(analytic functional) on $\Xi $. The converse statement is also
true: for any hyperfunction $\mu $ on $\Xi $, the function $u(x)$
in (\ref{E:represent}) is a solution of the equation $Pu=0$ in
$\Real^n$ which satisfies the growth condition
(\ref{E:hyper_estimate})
\end{enumerate}
\end{theorem}

\proof The second statement of the theorem is proven in our paper
\cite{KuPin} with Y.~Pinchover. So, we concentrate now on the
proof of the first one. The proof consists of three major parts:
defining appropriate function spaces and interpreting the
corresponding class of solutions as functionals; proving
Paley-Wiener type theorems for this class of spaces (Lemma
\ref{L:P-W} below); constructing a specific exact sequence of
topological spaces. The last step, i.e. constructing and proving
exactness of a sequence (Lemma \ref{L:isomorphism}) is usually the
most technical one.

Let us make first of all the following remark:
\begin{remark}
\label{R:pointwiseL2} {\em Using a standard elliptic argument
(Schauder type estimate) and periodicity of the equation, it is
standard to show that a solution satisfies for some $N$ the
pointwise growth condition (\ref{E:distr_estimate}) if and only if
it satisfies for some (different) $N$ the following $L_2$ growth
condition:
\begin{equation}\label{E:L2estimate}
u(x)(1+|x|)^{-N}\mathrm{e}^{ -h(x/\left| x\right| ) \left|
x\right|} \in L^2(\Bbb{R}^n)\,.
\end{equation}}
\end{remark}

Let us now return to the proof of the theorem. Assume first that a
function $u$ has the representation (\ref {E:represent}) with a
distribution $\mu$. Then it is obvious that it is a solution of
the equation $Pu=0$. We only need to establish the estimate
(\ref{E:distr_estimate}). Due to compactness of $\Xi$, the
distribution $\mu$ can be represented as a finite sum of terms of
the form $D^k\mu_k(\xi)$, where $D^k$ is a constant coefficient
homogeneous linear differential operator of order $k$ with respect
to $\xi$ and $\mu_k$ is a measure on $\Xi$. So, it is sufficient
to establish (\ref{E:distr_estimate}) for such a term only. In
other words, we need an estimate of the function $v(x)=\langle
\mu_k(\xi),D^k_\xi u_\xi(x)\rangle$. According to Lemma
\ref{L:Floq_struct}, $D^k_\xi u_\xi(x)$ is an analytically
depending on $\xi \in \Xi$ Floquet solution of $Pu=0$ of order
$k$. This means that it satisfies an estimate of the type
(\ref{E:distr_estimate}) with $N=k$. Then the estimate for $v(x)$
follows, since $\mu_k$ is a finite measure. Hence $u(x)$, being
the sum of a finitely many such terms, also satisfies
(\ref{E:distr_estimate}).

Suppose now that $u$ satisfies (\ref{E:distr_estimate}). We need
to prove that $u$ can be represented as in (\ref{E:represent}). In
order to do so, we need first to interpret this class of solutions
in dual terms.

Consider the following Fr\'{e}chet spaces of test functions:
$$ W_{m}:=\left\{ \phi \in H_{\mathrm{loc}}^m(\Real^n)\,|\,
<\phi>_{m,N}\, <\infty  \quad\forall N>0 \right\},
$$
where
$$<\phi>_{m,N}\,: =\sup_{\gamma \in \Gamma }\left\{\left|
\left| \phi \right|\right|_{H^m(K+\gamma)}(1+|\gamma|)^N
\mathrm{e}^{h(\gamma /\left| \gamma \right| )\left| \gamma
\right|}\right\} \, .
$$

The operator $P^{*}$ clearly maps continuously $W_{2} $ into
$W_{0}$. It is also clear that due to (\ref{E:distr_estimate}),
the linear functional
$$
<u,\phi >:=\int_{\mathbb{R}^n} u(x)\phi (x)\,\mathrm{d}x
$$
is continuous on the space $W_{0}$. Since $Pu=0$, Schauder
elliptic estimates together with the periodicity of the operator
show that estimates similar to (\ref{E:distr_estimate}) hold also
for the derivatives of $u$. One observes by a simple argument that
$u$ is a continuous functional on $W_{0}$, which annihilates the
range of the dual operator $P^{*}:W_{2}\rightarrow W_{0}$. Now we
can apply Floquet theory arguments analogous to the ones used in
\cite[Section 3.2]{Ku} or in \cite{KuPin} to obtain
(\ref{E:represent}). However, some technical details needed in the
cases considered in \cite{Ku,KuPin} and in this paper are
significantly different, so we provide the details of this
derivation.

First of all, we need to obtain a Paley-Wiener type theorem for
the Floquet transform in the spaces $W_{m}$. Let us denote by
$V^*$ the tube that consists of all points $z\in (\CC^*)^n$ such
that $z^{-1}=(z_1^{-1},...,z_n^{-1})\in V$, where the tube $V$ is
defined in (\ref{tubulard}). We introduce the space
$A^\infty(V^*)$ of holomorphic functions on the tube $V^*$ that
are infinitely differentiable up to its boundary $\partial V^*$.
Analogously, if $\mathcal{E}$ is a holomorphic Banach bundle in a
neighborhood of $V^*$, we denote by $A^\infty(V^*,\mathcal{E})$
the space of sections of $\mathcal{E}$ over the (closed) tube
$V^*$ that are holomorphic in the interior and infinitely
differentiable up to the boundary of $V^*$. This space is equipped
with the natural Fr\'{e}chet space topology. The following
statement is a Paley-Wiener type theorem for the transform
${\mathcal U}$ in the spaces $W_m$.
\begin{lemma}\label{L:P-W}
\begin{enumerate}
\item  The operator
$$
{\mathcal U}:W_{m}\rightarrow A^\infty (V^{*},{\mathcal E}_m).
$$
is a topological isomorphism.

\item Under the transform ${\mathcal U}$, the operator
$$
P^{*}:W_{2}\rightarrow W_{0}
$$
becomes the operator ${\mathcal P}(z)$ of multiplication by a
holomorphic Fredholm morphism between the fiber bundles ${\mathcal
E}_2$ and ${\mathcal E}_0$:
$$
A^\infty (V^{*},{\mathcal E}_2)\stackrel{{\mathcal P}(z)}{\rightarrow }%
A^\infty  (V^{*},{\mathcal E}_0). $$ Here ${\mathcal P}(z)$ acts
on each fiber of ${\mathcal E}_2$ as the restriction to this fiber
of the operator $ P^{*}$ acting between $H^2(K)$ and $L^2(K)$.
\end{enumerate}
\end{lemma}
Before proving this lemma, we first obtain the following auxiliary
statement:
\begin{lemma}\label{L:P-W-Hilbert}
Let $H$ be a complex Hilbert space and $W(H)$ be the Fr\'{e}chet
space of sequences $f=\{f_\g\},f_\g \in H, \g \in \G$ such that
the semi-norm
$$ \phi_{N}(f):= \sup_{\g \in \G}\left\{\|f_\g\|_{H}(1+|\gamma|)^N
\mathrm{e}^{h(\gamma /\left| \gamma \right| )\left| \gamma
\right|}\right\}
$$
is finite for any $N$. Here, as before, $h$ is the indicator
function (\ref{indicator}).

Then a sequence $f=\{f_\g\}$ belongs to $W(H)$ if and only if the
function
\begin{equation}\label{E:P-W-Hilbert}
\widehat{f}(z):=\sum\limits_{\g \in \G}f_{-\g} z^\g
\end{equation}
belongs to $A^\infty(V^*,H)$. The mapping $f \mapsto \widehat{f}$
is an isomorphism of the space $W(H)$ onto $A^\infty(V^*,H)$.
\end{lemma}
\pf Let $f \in W(H)$. We will show that the series
(\ref{E:P-W-Hilbert}) converges uniformly on $V^*$ as a series of
$H$-valued functions on $V^*$. This will imply that $\widehat{f}$
is analytic in $V^*$ and continuous up to the boundary. Then we
will check that the same holds for the derivatives of the series,
which will imply that $\widehat{f}\in A^\infty(V^*,H)$.

Taking into account that any $z\in V^*$ can be represented as
$z=\mathrm{e}^{-\mathrm{i}k}$ with $\mathrm{Im}k\in G$, and thus
$\mathrm{Im} k\cdot
   \g\leq h(\g/|\g|)|\g|$, we can estimate
\begin{equation}\label{E:Hilbert-estimate}
\begin{array}{cc}
   \|\widehat{f}(z)\|\leq\sum\limits_{\g \in \G}\|f_{-\g}\|\mathrm{e}^{-\mathrm{Im} k\cdot
   \g}
   =\sum\limits_{\g \in \G}\|f_{\g}\|\mathrm{e}^{\mathrm{Im} k\cdot
   \g}\\[4mm]
  \leq\sum\limits_{\g \in \G}(1+|\g|)^{-n-1}\|f_{\g}\|(1+|\g|)^{n+1}\mathrm{e}^{h(\g/|\g|)|\g|} \\[2mm]
   \leq \left[\sum\limits_{\g \in
   \G}(1+|\g|)^{-n-1}\right]\phi_{n+1}(f).
\end{array}
\end{equation}
Since the series $\sum_{\g \in \G}(1+|\g|)^{-n-1}$ converges, this
implies the analyticity in $V^*$ and continuity up to the boundary
of $\widehat{f}(z)$. Multiple differentiation with respect to $k$
amounts to multiplying the coefficients of (\ref{E:P-W-Hilbert})
by a polynomial with respect to $\g$ factor. Due to the definition
of the space $W(H)$, one can get an estimate from above similar to
(\ref{E:Hilbert-estimate}), but with the seminorm
$\phi_{n+d+1}(f)$ instead of $\phi_{n+1}(f)$, where $d$ is the
order of differentiation. Thus, in fact the function is infinitely
smooth up to the boundary. These estimates also prove that the
mapping $f\in W(H) \mapsto \widehat{f} \in A^\infty(V^*,H)$ is
continuous.

Let us now prove the surjectivity of this mapping. Assume that
$s(z)\in A^\infty (V^*,H)$. Let $z=\exp \mathrm{i}k$, then $s$ as
a function of $k$ is periodic with respect to the reciprocal
lattice $\Gamma ^{*}$. Expanding it into the Fourier series, we
get
\begin{equation}\label{E:Hilb_Fourier_expansion}
s(z)=\sum\limits_{\gamma \in \Gamma }s_{-\gamma }z^\gamma=\sum\limits_{\gamma \in \Gamma }s_{\gamma }z^{-\gamma},
\end{equation}
where $s_{-\gamma }\in H$. We need to show now that $\{s_\g\}\in W(H)$.
For this purpose, we use
the standard formulas for the Fourier coefficients:
$$
s_{\gamma }=\frac 1{(2\pi )^n}\int\limits_B
s(\mathrm{e}^{\mathrm{i}(\gb-\mathrm{i}\alpha
)})\mathrm{e}^{\mathrm{i}(\gb-\mathrm{i}\alpha )\cdot \gamma
}\,\mathrm{d}\gb,\quad\forall\, \alpha\in G,
$$
where $B$ is the first Brillouin zone, and we write $z=\exp
\mathrm{i}k=\exp [\mathrm{i}(\gb -\mathrm{i}\ga)], \ga\in G$.

Integrating by parts $l$ times with respect to $\gb$, where $l=(l_1,...,l_n)$ is a multi-index, we obtain analogously
\begin{equation}\label{E:Fourier_Hilbert}
s_{\gamma }=\frac {(-\mathrm{i}\g)^{-l}}{(2\pi )^n}\int\limits_B
\frac{\partial^l s}{\partial \gb
^l}(\mathrm{e}^{\mathrm{i}(\gb-\mathrm{i}\alpha
)})\mathrm{e}^{\mathrm{i}(\gb-\mathrm{i}\alpha )\cdot \gamma
}\,\mathrm{d}\gb\qquad\forall \alpha\in G.
\end{equation}
Now straightforward norm estimate in (\ref{E:Fourier_Hilbert}) gives
\begin{equation}\label{sgest}
\|s_\g\|_H \leq C \max_{z\in V^*}\| \frac{\partial^l s}{\partial
\gb ^l}(z)\| _{H}\,\g^{-l}\mathrm{e}^{-\ga\cdot\gamma}
\end{equation}
for any multi-index $l$ and any $\alpha\in
G$. Optimizing with respect to $\alpha\in G$, we get
\begin{equation}\label{sgest1}
\|s_\g\|_H \leq C_N (1+|\g|)^{-N}\mathrm{e}^{-h(\gamma /| \gamma |
)|\g |}
\end{equation}
for any $N$. This means that $f:=\{s_\g\}$ belongs to $W(H)$ and by its construction $\widehat{f}=s(z)$.
This proves Lemma \ref{L:P-W-Hilbert}. \qed

Let us now complete the proof of Lemma \ref{L:P-W}.

We start proving the first claim of the lemma. Let a function
$F(x)$ belong to $W_m$. Consider a sequence $f=\{f_\g\}$ of
elements of $H^m(K)$ defined as follows:
$$
f_\g(x)=F(x+\g)\qquad x\in K, \g\in \G.
$$
Then clearly the condition $F\in W_m$ is equivalent to two
conditions: the first one that $f\in W(H^m(K))$, and second that
$F\in H^m_{\mathrm{loc}}(\RR^n)$, i.e. that the functions $f_\g$
defined on shifted copies of the fundamental domain $K$, fit
smoothly across the boundaries.

Analogously, the requirement that a section $\phi$ belongs to
$A^\infty (V^*,{\mathcal E}_m)$ consists of two conditions. The
first one that $\phi\in A^\infty (V^*,H^m(K))$ and the second that
it is a section of the subbundle ${\mathcal E}_m\subset V^*\times
H^m(K)$.

We can notice now that the Floquet transform on $W_m$ is the
restriction of the transform $f\mapsto \widehat{f}$ of Lemma
\ref{L:P-W-Hilbert} from the larger space $W(H^m(K))$. Thus, Lemma
\ref{L:P-W-Hilbert} claims that this transform is an isomorphism
of $W(H^m(K))$ onto $A^\infty (V^*,H^m(K))$. On the other hand,
the second conditions: the fitting of $f_\g$ across the boundaries
and being a section of the subbundle ${\mathcal E}$, are
intertwined by the Floquet transform, according to the first
statement of Lemma \ref{Plancherel}. This proves the first claim
of the lemma.

Now, the second claim of Lemma \ref{L:P-W} follows from the third
one of Lemma \ref{Plancherel}. Lemma \ref{L:P-W} is proven. \qed

Let us now return to the proof of Theorem \ref{T:distribution}. We
remind the reader that we have a solution $u$ with the estimate
(\ref{E:distr_estimate}), for which we need to prove the
representation (\ref{E:represent}). Let us apply the Floquet
transform ${\mathcal U}$. Then the image ${\mathcal U}u$ of the
solution $u$ under the Floquet transform is a continuous linear
functional on $A^\infty (V^{*},{\mathcal E}_0)$, which is in the
cokernel of the operator
$$
A^\infty  (V^{*},{\mathcal E}_2)\stackrel{{\mathcal P}(z)}{\rightarrow }%
A^\infty  (V^{*},{\mathcal E}_0).
$$
This, indeed is a one-to-one correspondence between solutions of
the required class and such functionals. Thus, we need to describe
all such functionals. Let $u_z(\cdot)=z^xp(z,\cdot)$ be the Bloch
solution of the equation $Pu=0$ introduced in Lemma
\ref{Lambda-lemma2}. We will also employ the space $C^\infty
(\Psi)$ with the standard topology, where the smooth variety
$\Psi$ is introduced in (\ref{Psi}). Consider the mapping
$$
t:A^\infty (V^{*},{\mathcal E}_0)\rightarrow C^\infty (\Psi)
$$
that for a section $f(z,x) \in A^\infty (V^{*},{\mathcal E}_0)$ of
the bundle ${\mathcal E}_0$ produces
$$
t_{f}(z)=<f(z^{-1},\cdot),u_{z}(\cdot)>=\int\limits_{\mathbb{T}^n}f(z^{-1},x)u_{z}(x)\,\mathrm{d}x.
$$
Here $z^{-1}=(z_1^{-1},\ldots, z_n^{-1})$.

As we will see soon, the following lemma will finish the proof of
the theorem:
\begin{lemma} \label{L:isomorphism} The
mapping $t$ is a topological homomorphism and the following
sequence is exact:
\begin{equation}
A^\infty (V^{*},{\mathcal E}_2)\stackrel{{\mathcal P}(z)}{\rightarrow }%
A^\infty (V^{*},{\mathcal E}_0)\stackrel{t}{\rightarrow }C^\infty
(\Psi) \rightarrow 0. \label{E:sequ}
\end{equation}
\end{lemma}
{\bf Proof of the lemma.} Continuity of ${\mathcal P}(z)$ is
already established. Continuity of $t$ is obvious. The complex
property of the sequence (\ref{E:sequ}) (i.e. that $t{\mathcal
P}(z)=0$) follows from the construction of $t$. Thus, the only
thing that requires proof is exactness in the second and third
terms. The topological homomorphism property will follow then from
exactness and the open mapping theorem. So, we only need to prove
that: i) any section $\phi\in A^\infty(V^*,{\mathcal E}_0)$ such
that $t\phi=0$ belongs to the range of ${\mathcal P}(z)$ and ii)
any function $f\in C^\infty(\Psi)$ is in the range of $t$.

Let us start with the first of these tasks. So, let $\phi\in
A^\infty(V^*,{\mathcal E}_0)$ be such that $t\phi =0$. Consider
the inverse ${\mathcal P}^{-1}(z)$ to the morphism ${\mathcal
P}(z)$. It is defined (and hence holomorphic) in a neighborhood
$V_\epsilon^*$ of the tube $V^*$, except for an analytic
submanifold, whose intersection with $V^*$ is $\Psi$ (see Lemma
\ref{Lambda-lemma2}). Let us consider the function $f={\mathcal
P}^{-1}(z)\phi(z)$. The only thing now to prove is that this
function does not have any singularities along $\Psi$. This is a
local question, so let us return in a neighborhood of a point of
$\Psi$ to the quasi-momenta coordinates $k$ and consider the
structure of the inverse ${\mathcal P}^{-1}(z)$. As it was shown
in the proof of \cite[Lemma 21]{KuPin}, the inverse has the form
$B(k)/\Lambda (k)$, where $B(k)$ is an analytic operator-valued
function. This means that $f(k)=(B(k)\phi (k))/\Lambda (k)$. The
condition $t\phi=0$ guarantees that the numerator $g(k)=B(k)\phi
(k)\in A^\infty(V^*,H)$ vanishes on $\Psi$, where $H$ is a Hilbert
space. Our goal is to prove that this is sufficient for its smooth
divisibility (on $\partial V^*$) by $\Lambda$. We recall here that
$\Lambda$ is analytic in a vicinity of $\partial V^*$ and has
simple zeros along $\Psi$ (Lemmas \ref{Lambda-lemma} and
\ref{Fermi-lemma}). We notice that it is sufficient to prove this
for scalar functions, i.e. for $H=\CC$. This can be justified in
many different ways. For instance, the statement is local, and
locally, due to the Fredholm nature of the morphism ${\mathcal
P}(z)$, one can project the problem onto a finite dimensional
subspace, using a lemma by M.~Atiyah \cite{Atiyah} (see also
\cite[Lemma 2.1]{ZK} and \cite[Lemma 1.2.11 and Theorem
1.3.9]{Ku}), which will reduce it to a finite dimensional, and
thus also to scalar case. So, we will assume in this part of the
proof that $g\in A^\infty$ is a scalar function. According to a
result of \cite{Horm,Lo} (see also \cite{Malgr2,Pa} and
\cite[Theorem 1.1' in Ch. VI]{Malgr2}), it is sufficient to check
the divisibility at each point of $\Psi$ on the level of formal
Taylor series. So, let us pick a point $k$ of $\Psi$ and introduce
coordinates $x\in \RR^{n-1}$ in the tangent space $T_k (\Psi)\in
i\RR^n$. The complexification $T^c_k(\Psi)$ of this tangent space
is a part of the tangent space to the boundary of the tube. Let us
chose coordinates $y\in \RR^{n-1}$ in $T^c_k (\Psi) \cap \RR^n$
that correspond to the coordinates $x$ in $T_k (\Psi)$. An extra
coordinate $t$ in $T_k (\Psi) \cap \RR^n$ is required to obtain
the whole tangent space $T_k (\partial V^*)$. Let us denote by
$\widehat{g}(x,y,t)$ and $\widehat{\Lambda}(x,y,t)$ the formal
Taylor series of $g$ and $\Lambda$ at the point $k$. Then we know
that $\widehat{g}(x,0,0)=0$ and $\widehat{\Lambda}(x,0,0)=0$
(formal power series versions of vanishing of functions $g$ and
$\Lambda$ on $\Psi$). Recall that $\widehat{g}(x,y,t)$ is the
series for a CR-function $g$ on the boundary (since $g$ is the
boundary value of an analytic function). This means that
$\widehat{g}(x,y,t)$ satisfies Cauchy-Riemann conditions with
respect to the variable $z=x+iy\in \CC^{n-1}$. Then uniqueness of
analytic continuation\footnote{The uniqueness of analytic
continuation in this power series setting is straightforward to
derive algebraically directly from the Cauchy-Riemann conditions
for power series.} claims that $\widehat{g}(x,0,0)=0$ for all $x$
implies $\widehat{g}(x,y,0)=0$ for all $(x,y)$. The same is true
for $\widehat{\Lambda}$, due to analyticity of $\Lambda$. Now, in
coordinates $z=x+iy,t$ we are dealing with the formal series
$\widehat{g}(z,t)$ and $\widehat{\Lambda}(z,t)$, both of which
vanish at $t=0$ and such that $\widehat{\Lambda}$ has zero of
first order at $t=0$. Then, vanishing of $\widehat{g}(z,0)$
guarantees divisibility in formal series of $\widehat{g}$ by
$\widehat{\Lambda}$. As it was explained above, this implies
smooth divisibility of $g$ by $\Lambda$ and thus finishes the
proof of exactness in the second term of the sequence
(\ref{E:sequ}).

Let us now prove the exactness in the third term of the sequence.
First of all, we notice that the vector-function $u_z$, defined on
$\Psi$ only, can be extended to an analytic vector-function (which
we will denote the same way) on $V^*_\epsilon$ for some small
epsilon. Indeed, as it is shown in \cite{KuPin}, $V^*_\epsilon$ is
a Stein manifold. Then, according to the Corollary 1 from the
Bishop's theorem 3.3 in \cite{ZK} (see the original theorem in
\cite{Bi}), the restriction mapping to an analytic subset of a
Stein variety is surjective. Thus, the required extension of $u_z$
exists. Let also $v(z)$ be a holomorphic family such that
$tv(z)|_\Psi=1$ (it is not hard to prove the existence of such a
family). Consider a function $\phi(z)\in C^\infty (\Psi)$. Notice
that the domain $V^*$ is strictly pseudo-convex and the
complexifications of the tangent spaces to the submanifold
$\Psi\subset \partial V^*$ are parts of the tangent spaces to
$\partial V^*$. Thus, $\Psi$ and $\partial V^*$ satisfy the
conditions of \cite{ChaCho2} needed for $\Psi$ to be an $A^\infty$
interpolation variety, and hence the restriction mapping
$A^\infty(V^*)\mapsto C^\infty (\Psi)$ is surjective. Hence, there
exists a function $\psi\in A^\infty(V^*)$ such that
$\psi|_\Psi=\phi$. Now taking $f=\psi (z) v(z)\in
A^\infty(V^*,{\mathcal E}_0)$ guarantees that $tf=\phi$. This
finishes the proof of the lemma. \qed

It is easy now to finish the proof of the theorem. Indeed, after
the Floquet transform solution $u$ becomes a continuous linear
functional on $A^\infty (V^{*},{\mathcal E}_0)$ that annihilates
the range of the operator of multiplication by ${\mathcal P}(z)$.
Lemma \ref{L:isomorphism} implies that such a functional can be
pushed down to the space $C^\infty (\Psi)$. Any such functional is
a distribution $\mu $. Hence, the action $<u,\phi >$ of the
functional $u$ on a function $\phi\in W_{0}$ can be obtained as
$$
<u,\phi >=<\mu (z),t(z)({\mathcal U}\phi )>.
$$
Applying now the explicit formulas for the transforms ${\mathcal
U}$ and $t$, one arrives to the representation
(\ref{E:represent}). Indeed,
\begin{eqnarray} \label{range}
t_{({\mathcal U}\phi )}(z)=\int\limits_{K}{\mathcal U}\phi
(z^{-1},x)u_{z}(x)\,\mathrm{d}x  \\ =\sum \limits_{\gamma \in
\Gamma}\int\limits_{K-\gamma } \phi
(x)z^{-\gamma}u_{z}(x+\gamma)\,\mathrm{d}x \nonumber \\
=\int\limits_{\Real^n }\phi (x)u_{z}(x)\,\mathrm{d}x. \nonumber
\end{eqnarray}
In this calculation we used the property of the Bloch solutions $$
u_{z}(x+\gamma)=z^{\gamma}u_{z}(x). $$ Therefore, $$
<u,\phi >=<<\mu (z),u_{z}>,\phi >,
$$
which concludes the proof of the theorem. \qed

\begin{center}
{\bf Acknowledgments} \end{center} The author expresses his
gratitude to Y.~Pinchover, the co-author of the previous papers
\cite{KuPin,KuPin2}, with whom this manuscript has been discussed
on numerous occasions and who has made many suggestions that have
improved the text, and to A.~Tumanov for helpful information on
pick sets results.

The work of the author was partially supported by the NSF Grant
DMS 0406022 and by Grant No. 1999208 from the United States-Israel
Binational Science Foundation (BSF). The author expresses his
gratitude to the NSF and BSF for this support. The content of this
paper does not necessarily reflect the position or the policy of
the federal government of the USA, and no official endorsement
should be inferred.

\end{document}